\documentclass[3p,times]{elsarticle}

\usepackage{amsfonts}
\usepackage{mathrsfs}
\usepackage{amsmath}
\usepackage{amssymb}

 \usepackage{lineno}
\newtheorem{theorem}{\indent Theorem}[section]
\newtheorem{corollary}{\indent Corollary}[section]

\newtheorem{lemma}{\indent Lemma}[section]

\journal{Discrete and Continuous Dynamical Systems-B}

\begin{document}

\begin{frontmatter}



\title{Long-time behavior of the three dimensional globally modified Navier-Stokes equations}

\author{Fang Li$^a,$ Bo You$^{b,*}$}

\address{$^a$ School of Mathematics and Statistics, Xidian University, Xi'an, 710126, P. R. China\\
$^b$ School of Mathematics and Statistics, Xi'an Jiaotong University, Xi'an, 710049, P. R. China}

 \cortext[cor1]{Corresponding author.\\
 E-mail addresses:fli@xidian.edu.cn (F. Li),youb2013@xjtu.edu.cn(B. You).}


\begin{abstract}
This paper is concerned with the long-time behavior of solutions for the three dimensional globally modified Navier-Stokes equations in a
three-dimensional bounded domain. We prove the existence of a global attractor $\mathcal{A}_0$ in $H$ and investigate the regularity of the global attractors by proving that $\mathcal{A}_0=\mathcal{A}$ established in \cite{ct}, which implies the asymptotic smoothing effect of solutions for the three dimensional globally modified Navier-Stokes equations in the sense that the solutions will eventually become more regular than the initial data. Furthermore, we construct an exponential attractor in $H$ by verifying the smooth property of the difference of two solutions, which entails the fractal dimension of the global attractor is finite.
\end{abstract}

\end{frontmatter}


\section{Introduction}
\def\theequation{1.\arabic{equation}}\makeatother
\setcounter{equation}{0}
In this paper, we consider the long-time behavior of solutions for the following three dimensional globally modified Navier-Stokes equations (see \cite{ct})
\begin{equation}\label{1.1}
\begin{cases}
&\frac{\partial u}{\partial t}-\nu\Delta u+F_N(\|u\|)(u\cdot\nabla)u+\nabla p=f,\,\,(x,t)\in\Omega\times\mathbb{R}^+,\\
 &\nabla\cdot u=0,\quad\quad\quad\,\,\quad\quad\quad\quad\quad\quad\quad\quad\,\,(x,t)\in\Omega\times\mathbb{R}^+,\\
 & u=0,\quad\quad\quad\,\,\quad\quad\quad\,\,\quad\quad\quad\quad\quad\quad\,\,\,(x,t)\in\partial\Omega\times\mathbb{R}^+,\\
 &u(x,0)=u_0(x),\quad\quad\quad\,\,\quad\quad\quad\,\,\quad\quad\,\,\,\,\,x\in\Omega,
\end{cases}
\end{equation}
where $\Omega\subset\mathbb{R}^3$ is a bounded domain with smooth boundary $\partial\Omega$ and $\mathbb{R}^+=[0,+\infty),$ $F_N(r)$ is defined by
\begin{align*}
F_N(r)=\min\{1,\frac{N}{r}\},~r\in\mathbb{ R}^+
\end{align*}
for some $N \in \mathbb{R}^+$ and $\|u\|$ is the norm of $u$ in $V$ defined in the next section, $u(x,t)$ is the velocity of the fluid and $p(x,t)$ is the pressure of the fluid, $\nu$ is the kinematic viscosity of the fluid,
$f$ is the given external force.

 In the past several decades, the well-posedness and the long-time behavior of solutions for the three dimensional globally modified Navier-Stokes equations have been extensively studied from the theoretical point of view (see \cite{ct3, ct, ct1, ct2, dbq, kpe, kpe1, mp, mp2, mp3, mp1, rm}). In particular, in \cite{ct}, the authors have considered the existence and uniqueness of strong solutions of the three dimensional globally modified Navier-Stokes equations. Meanwhile, they also established the existence of a global attractors in $V$ by verifying the flattening property and obtained the existence of bounded entire weak solutions of the three dimensional Navier-Stokes equations with time independent forcing via a limiting argument. The uniqueness of weak solutions for the three dimensional globally modified Navier-Stokes equations was proved in \cite{rm}. In \cite{mp3}, the authors have established the existence, uniqueness and continuity properties of solutions for the three dimensional globally modified Navier-Stokes equations with finite delay terms within a locally Lipschitz operator. Moreover, they also proved that there exists a unique globally asymptotically exponentially stable stationary solution for the stationary problem under suitable assumptions. The regularity of solutions as well as the relationship between global attractors, invariant measures, time-average measures and statistical solutions of the three dimensional globally modified Navier-Stokes equations were investigated in the case of temporally independent forcing in \cite{ct1}. In \cite{kpe}, the authors have established the existence and finite fractal dimension of a pullback attractor in $V$ for the three-dimensional non-autonomous globally modified Navier-Stokes equations in a bounded domain under some appropriate assumptions on the time dependent forcing term. The existence and uniqueness of strong solutions, asymptotic behavior of solutions and the existence of a pullback attractor for the three-dimensional globally modified Navier-Stokes equations with delay in the locally Lipschitz case was considered in \cite{ct2}. In \cite{mp}, the authors have established the existence and uniqueness of solution for the three dimensional globally modified Navier-Stokes equations containing infinite delay terms. Moreover, they also proved the global exponential decay of the solutions of the evolutionary problem to the stationary solution under some suitable additional conditions. The existence of pullback attractors in two different settings for the three dimensional globally modified Navier-Stokes equations with infinite delays and their relationship were investigated in \cite{mp2}. In \cite{kpe1}, the authors have established the relationship between invariant measures and statistical solutions for the three-dimensional non-autonomous globally modified Navier-Stokes equations in the case of temporally independent forcing term by using a smoother Galerkin scheme and proved that a stationary statistical solution is also an invariant probability measure under suitable assumptions.

 Although the global attractor represents the first important step
in the understanding of long-time behavior of dynamical systems generated by the three dimensional globally modified Navier-Stokes equations, it may also present two essential drawbacks:on the one hand, the rate of attraction of the trajectories may be small and it is usually very difficult to estimate this rate in terms of the physical parameters of the problem. On the other hand, it is very sensitive to perturbations such that the global attractor can change drastically under very small perturbations of the initial dynamical system. These drawbacks obviously lead to essential difficulties in numerical simulations of global attractors and even make the global attractor unobservable in some sense.

An alternative object to describe the long-term dynamics is an inertial manifold, which is free from the above-mentioned drawbacks (see \cite{em1}). Unfortunately, its existence can be proved only under very restrictive spectral gap assumptions, which can be verified in few particular dynamical systems, mainly arising from one-dimensional parabolic equations (see \cite{ha}).

In order to overcome this difficulty, an intermediate object has been introduced in \cite{cvv, ea}, an exponential attractor or inertial set. The exponential attractors contain the global attractor, are finite dimensional, and attract the trajectories exponentially fast. In contrast to the global attractor, an exponential attractor attracts exponentially the trajectories and will thus be more stable. Meanwhile, it also provides a way of proving that the global attractor has finite fractal dimension. Furthermore, in some situations, the global attractor can be very simple and thus fails to capture interesting transient behaviors. In such situations, an exponential attractor could be a more suitable object. Therefore, it is useful to explore the existence of an exponential attractor for the three dimensional globally modified Navier-Stokes equations.

The purpose of this paper is to study long-time behavior of solutions for the three dimensional globally modified Navier-Stokes equations. In the next section, we give the definition of some function spaces and recall a lemma used in the sequel. In Section 3, we first prove the existence of a global attractor and the regularity of global attractors, and then, we construct an exponential attractor for problem \eqref{1.1} by proving the smooth property of the difference of two solutions, which implies the fractal dimension of the global attractor is finite.

 Throughout this paper, let $X$ be a Banach space endowed with the norm $\|\cdot\|_X,$ let $\|u\|_p$ be the norm of $u$ in $L^p(\Omega)$ and let $C$ be positive constants which may be different from line to line.
\section{Preliminaries}
\def\theequation{2.\arabic{equation}}\makeatother
\setcounter{equation}{0}
In this section, we introduce some function spaces and recall a lemma used in the sequel.

Define
\begin{align*}
\mathcal {V}=\{v\in (C_0^{\infty}(\Omega))^3:\nabla\cdot v=0\}.
\end{align*}
Denote by $H$ and $V,$ respectively, the closure of $\mathcal{V}$ with respect to the $(L^2(\Omega))^3$-norm $\|\cdot\|_2$ and the $(H_0^1(\Omega))^3$-norm $\|\cdot\|.$

\begin{lemma}(\cite{rm})\label{2.1}
For every $u,$ $v\in V$ and each $N>0,$ the following two conclusions hold:
\begin{itemize}
\item [(i)] $0\leq \|u\| F_N(\|u\|)\leq N.$
\item [(ii)] $|F_N(\|u\|)-F_N(\|v\|)|\leq \frac{1}{N}F_N(\|u\|)F_N(\|v\|)\|u-v\|.$
\end{itemize}
\end{lemma}
\section{The existence of a global attractor and an exponential attractor}
\def\theequation{3.\arabic{equation}}\makeatother
\setcounter{equation}{0}
\subsection{The well-posedness of weak solutions}
The existence and uniqueness of weak solutions was obtained in \cite{ct, rm}. Here, we only state it as follows.
\begin{theorem}\label{3.1.1}
Assume that $f\in H.$ Then for any $u_0\in H,$ there exists a unique weak solution $u\in C(\mathbb{R}^+;H)\cap L^2_{loc}(\mathbb{R}^+;V)$ of problem \eqref{1.1}, which depends continuously on the initial data in $H.$
\end{theorem}

   By Theorem \ref{3.1.1}, we can define the operator semigroup $\{S(t)\}_{t\geq 0}$ in $H$ by
\begin{align*}
S(t)u_0:=u(t)=u(t,u_0)
\end{align*}
for any $t\geq0,$ where $u(t)$ is the weak solution of problem \eqref{1.1} with initial data $u(0)=u_0 \in H.$
\subsection{The existence of a global attractor}
In this subsection, we will prove the existence of a global attractor in $H$ and the regularity of global attractors for problem \eqref{1.1}. In what follows, we first prove the existence of an absorbing set in $V$ for problem \eqref{1.1}.
\begin{theorem}\label{3.2.1}
Assume that $f\in H.$ Then there exists a positive constant $\rho_1$ satisfying for any bounded subset $B\subset H,$ there exists some time $t_0=t_0(B)>0$ such that for any weak solutions of problem \eqref{1.1} with initial data $u_0\in B,$ we have
\begin{align*}
\|u(t)\|^2\leq\rho_1
\end{align*}
for any $t\geq t_0.$
\end{theorem}
\textbf{Proof.} Taking the inner product of the first equation of \eqref{1.1} with $u$ in $H,$ we obtain
\begin{align*}
\frac{1}{2}\frac{d}{dt}\|u(t)\|^2_2+\nu\|u(t)\|^2=&\int_{\Omega}f\cdot u(t)\,dx\\
\leq&\|f\|_2\|u(t)\|_2\\
\leq&\frac{\nu}{2}\|u(t)\|^2+\frac{1}{2\nu\lambda_1}\|f\|_2^2.
\end{align*}
From the Poincar\'{e} inequality and Young inequality, we deduce
\begin{align}\label{3.2.2}
\frac{d}{dt}\|u(t)\|^2_2+\nu\|u(t)\|^2\leq\frac{1}{\nu\lambda_1}\|f\|_2^2
\end{align}
and
\begin{align*}
\frac{d}{dt}\|u(t)\|^2_2+\nu\lambda_1\|u(t)\|_2^2\leq\frac{1}{\nu\lambda_1}\|f\|_2^2,
\end{align*}
where $\lambda_1$ is the first eigenvalue of the Stokes operator $A.$

It follows from the classical Gronwall inequality that
\begin{align*}
\|u(t)\|_2^2\leq\|u_0\|_2^2 e^{-\nu\lambda_1 t}+\frac{1}{\nu^2\lambda_1^2}\|f\|_2^2,
\end{align*}
which implies that for any bounded subset $B\subset H,$ there exists some time $\tau_0=\tau_0(B)\geq 0$ such that for any $u_0\in B,$ we obtain
\begin{align}\label{3.2.3}
\|u(t)\|_2^2\leq \frac{2}{\nu^2\lambda_1^2}\|f\|_2^2
\end{align}
for any $t\geq\tau_0.$

Integrating \eqref{3.2.2} from $t$ to $t+1$ and using \eqref{3.2.3}, we obtain for any $u_0\in B,$
\begin{align}\label{3.2.4}
\nonumber\int_t^{t+1}\|u(s)\|^2\,ds\leq &\frac{1}{\nu^2\lambda_1}\|f\|_2^2+\frac{1}{\nu}\|u(t)\|_2^2\\
\leq&(1+\frac{2}{\nu\lambda_1})\frac{1}{\nu^2\lambda_1}\|f\|_2^2
\end{align}
for any $t\geq\tau_0.$

Multiplying the first equation of \eqref{1.1} by $Au=-P\Delta u$ and integrating over $\Omega,$ where $P$ is the Leray-Helmotz projection from $(L^2(\Omega))^3$ onto $H,$ we obtain
\begin{align*}
\frac{1}{2}\frac{d}{dt}\|u(t)\|^2+\nu\|Au(t)\|^2=&\int_\Omega f\cdot Au-F_N(\|u\|)[(u\cdot\nabla)u]\cdot Au\,dx\\
\leq&\|f\|_2\|Au\|_2+|F_N(\|u\|)|\,\|u\|_6\,\|\nabla u\|_3\,\|Au\|_2.
\end{align*}
It follows from Sobolev inequality, Young inequality and Lemma \ref{2.1} that
\begin{align*}
\frac{1}{2}\frac{d}{dt}\|u(t)\|^2+\nu\|Au(t)\|^2\leq&\|f\|_2\|Au\|_2+|F_N(\|u\|)|\,\|u\|_6\,\|\nabla u\|_3\,\|Au\|_2\\
\leq&\frac{1}{\nu}\|f\|_2^2+\frac{\nu}{4}\|Au\|_2^2+C|F_N(\|u\|)|\,\|u\|\,\|u\|^{\frac{1}{2}}\,\|Au\|_2^{\frac{3}{2}}\\
\leq&\frac{1}{\nu}\|f\|_2^2+\frac{\nu}{2}\|Au\|_2^2+CN^4\|u\|^2,
\end{align*}
which implies that
\begin{align}\label{3.2.5}
\frac{d}{dt}\|u(t)\|^2+\nu\|Au(t)\|^2\leq\frac{2}{\nu}\|f\|_2^2+CN^4\|u\|^2.
\end{align}
We deduce from the uniform Gronwall inequality and \eqref{3.2.4} that for any $u_0\in B,$
\begin{align}\label{3.2.6}
\nonumber\|u(t)\|^2\leq&\left(\frac{2}{\nu}\|f\|_2^2+(1+\frac{2}{\nu\lambda_1})\frac{1}{\nu^2\lambda_1}\|f\|_2^2\right)e^{CN^4}\\
\leq&\frac{1}{\nu}(2+\frac{1}{\nu\lambda_1}+\frac{2}{\nu^2\lambda_1^2})\|f\|_2^2 e^{CN^4}
\end{align}
for any $t\geq \tau_0+1.$

We infer from \eqref{3.2.5} that
\begin{align}\label{3.2.7}
\frac{d}{dt}\left((t-\tau_0)\|u(t)\|^2\right)\leq&\frac{2}{\nu}\|f\|_2^2(t-\tau_0)+CN^4(t-\tau_0)\|u\|^2+\|u(t)\|^2.
\end{align}
Integrating \eqref{3.2.7} from $\tau_0$ to $\tau_0+1$ and using \eqref{3.2.4}, we obtain for any $u_0\in B,$
\begin{align}\label{3.2.8}
\nonumber\|u(\tau_0+1)\|^2\leq&\int_{\tau_0}^{\tau_0+1}\frac{2}{\nu}\|f\|_2^2(t-\tau_0)+CN^4(t-\tau_0)\|u(t)\|^2+\|u(t)\|^2\,dt\\
\leq&\frac{2}{\nu}\|f\|_2^2+(CN^4+1)(1+\frac{2}{\nu\lambda_1})\frac{1}{\nu^2\lambda_1}\|f\|_2^2.
\end{align}
\qed\hfill

Thanks to the compactness of $V\subset H,$ we immediately obtain the following result.
\begin{theorem}\label{3.2.9}
Assume that $f\in H.$ Then the semigroup $\{S(t)\}_{t\geq 0}$ generated by problem \eqref{1.1} possesses a global attractor $\mathcal{A}_0$ in $H.$
\end{theorem}

Next, we will prove the regularity of global attractors for problem \eqref{1.1}.
\begin{theorem}\label{3.2.10}
Assume that $f\in H$ and let $\mathcal{A}$ be the global attractor in $V$ of problem \eqref{1.1} established in \cite{ct}. Then $\mathcal{A}=\mathcal{A}_0.$
\end{theorem}
\textbf{Proof.} From the definition of the global attractor, we deduce that $\mathcal{A}$ is a bounded subset of $V,$ which implies that $\mathcal{A}$ is also a bounded subset of $H.$ Therefore, we obtain
\begin{align*}
dist_H(\mathcal{A},\mathcal{A}_0)=&dist_H(S(t)\mathcal{A},\mathcal{A}_0)\,\,(\forall\,\,t\in\mathbb{R}^+)\\
=&\lim_{t\rightarrow+\infty}dist_H(S(t)\mathcal{A},\mathcal{A}_0)\\
=&0,
\end{align*}
which entails that
\begin{align}\label{3.2.11}
\mathcal{A}\subset\mathcal{A}_0.
\end{align}
Since $\mathcal{A}_0$ is a bounded subset of $H,$ we deduce from the proof of Theorem \ref{3.2.1} that there exists some time $t_1=t_1(\mathcal{A}_0)$ such that
\begin{align*}
S(t_1)\mathcal{A}_0\subset B_0\triangleq\left\{u\in V:\|u\|^2\leq\frac{2}{\nu}\|f\|_2^2+(CN^4+1)(1+\frac{2}{\nu\lambda_1})\frac{1}{\nu^2\lambda_1}\|f\|_2^2\right\},
\end{align*}
which implies that $S(t_1)\mathcal{A}_0$ is a bounded subset of $V.$ Hence, we have
\begin{align*}
dist_H(\mathcal{A}_0,\mathcal{A})=&dist_H(S(t)\mathcal{A}_0,\mathcal{A})\,\,(\forall\,\,t\in\mathbb{R}^+)\\
\leq&\frac{1}{\lambda_1}dist_V(S(t)\mathcal{A}_0,\mathcal{A})\,\,(\forall\,\,t\in\mathbb{R}^+)\\
=&\frac{1}{\lambda_1}\lim_{t\rightarrow+\infty}dist_V(S(t)\mathcal{A}_0,\mathcal{A})\\
=&\frac{1}{\lambda_1}\lim_{t\rightarrow+\infty}dist_V(S(t-t_1)S(t_1)\mathcal{A}_0,\mathcal{A})\\
=&0,
\end{align*}
entails that
\begin{align}\label{3.2.12}
\mathcal{A}_0\subset\mathcal{A}.
\end{align}
It follows from \eqref{3.2.11}-\eqref{3.2.12} that
\begin{align*}
\mathcal{A}_0=\mathcal{A}.
\end{align*}
\qed\hfill

\subsection{The existence of an exponential attractor}
In this subsection, inspired by the idea in \cite{ba, ea,em, rjc}, we prove the existence of an exponential attractor in $H$ for problem \eqref{1.1}. The definition of exponential attractor can be referred to \cite{ba, cvv, ea, em, em1, rjc}.

In what follows, we first prove the smoothing property of the semigroup $\{S(t)\}_{t\geq0}$ generated by problem \eqref{1.1}.
\begin{theorem}\label{3.3.1}
Assume that $f\in H$ and let $(u^i,p^i)$ be the solution of problem \eqref{1.1} with the initial data $u_0^i$ for $i=1,2.$ Then for any bounded subset $B\subset H,$ there exists some time $t_2=t_2(B)$ such that the following estimate holds
\begin{align}
\|u^1(t)-u^2(t)\|^2\leq \frac{\varrho_1}{\bar{t}}e^{\varrho_2 \bar{t}}\|u_0^1-u_0^2\|_2^2
\end{align}
for any $u_0^i\in B$ and any $t\geq t_2+1,$ where $\bar{t}=t-t_2,$ $\varrho_1$ and $\varrho_2$ are positive constants which only depend on $\Omega,$ $\nu,$ $\lambda_1,$ $N$ and $\|f\|_2.$
\end{theorem}
\textbf{Proof.} For any bounded subset $B$ of $H,$ we deduce from Theorem \ref{3.2.1} that there exists some time $t_2=t_2(B)$ such that
\begin{align*}
S(t)B\subset B_1\triangleq\{u\in V:\|u\|^2\leq \rho_1\}
\end{align*}
for any $t\geq t_2.$

For any $t\geq t_2,$ integrating \eqref{3.2.5} over $(t_2, t),$ we obtain that for any $u^1(0)\in B,$
\begin{align}\label{3.3.2}
\nonumber\nu\int_{t_2}^t\|Au^1(s)\|_2^2\,ds\leq&\frac{2}{\nu}\|f\|_2^2(t-t_2)+CN^4\int_{t_2}^t\|u^1(s)\|^2\,ds+\|u^1(t_2)\|^2\\
\leq&\frac{2}{\nu}\|f\|_2^2(t-t_2)+CN^4\rho_1(t-t_2)+\rho_1.
\end{align}

Let $u_0^i\in B$ for $i=1,2$ and $(u(t),p(t))=(u^1(t)-u^2(t),p^1(t)-p^2(t)),$ then $(u(t),p(t))$ satisfies the following equations
\begin{equation}\label{3.3.3}
\begin{cases}
&\frac{\partial u}{\partial t}-\nu\Delta u+F_N(\|u^1\|)(u^1\cdot\nabla)u^1-F_N(\|u^2\|)(u^2\cdot\nabla)u^2+\nabla p=0,\,\,(x,t)\in\Omega\times\mathbb{R}^+,\\
&\nabla\cdot u=0,\quad\quad\quad\,\,\quad\quad\quad\quad\quad\quad\quad\quad\quad\quad\quad\quad\quad\quad\quad\quad\quad\quad\,\,(x,t)\in\Omega\times\mathbb{R}^+,\\
&u=0,\quad\quad\quad\,\,\quad\quad\quad\quad\quad\quad\quad\quad\quad\quad\quad\quad\quad\quad\quad\quad\quad\quad\quad\,\,\,\,(x,t)\in\partial\Omega\times\mathbb{R}^+,\\
&u(x,0)=u_0^1-u_0^2,\quad\quad\quad\,\,\quad\quad\quad\quad\quad\quad\quad\quad\quad\quad\quad\quad\quad\quad\quad\,\,\,x\in\Omega.
\end{cases}
\end{equation}
Multiplying the first equation of \eqref{3.3.3} by $u$ and integrating over $\Omega,$ from the proof of Theorem 1.1 in \cite{rm}, we obtain
\begin{align*}
\frac{1}{2}\frac{d}{dt}\|u(t)\|_2^2+\nu\|u(t)\|^2\leq CN\|u(t)\|_2^{\frac{1}{2}}\|u(t)\|^{\frac{3}{2}}
\end{align*}
entails that
\begin{align*}
\frac{d}{dt}\|u(t)\|_2^2+\nu\|u(t)\|^2\leq CN^4\|u(t)\|_2^2.
\end{align*}
We infer from the classical Gronwall inequality that
\begin{align}\label{3.3.4}
\|u(t)\|_2^2+\nu\int_0^t\|u(s)\|^2\,ds\leq e^{CN^4t}\|u(0)\|_2^2.
\end{align}
Taking the inner product of the first equation of \eqref{3.3.3} with $Au$ in $H,$ we find
\begin{align}\label{3.3.5}
\frac{1}{2}\frac{d}{dt}\|u(t)\|^2+\nu\|Au(t)\|_2^2=-\int_\Omega\left(F_N(\|u^1\|)(u^1\cdot\nabla)u^1-F_N(\|u^2\|)(u^2\cdot\nabla)u^2\right)\cdot Au\,dx.
\end{align}
From the proof of Theorem 1.1 in \cite{rm}, we know
\begin{align*}
&F_N(\|u^1\|)(u^1\cdot\nabla)u^1-F_N(\|u^2\|)(u^2\cdot\nabla)u^2\\
=&F_N(\|u^1\|)(u\cdot\nabla)u^1+F_N(\|u^2\|)(u^2\cdot\nabla)u+\left(F_N(\|u^1\|)-F_N(\|u^2\|)\right)(u^2\cdot\nabla)u^1.
\end{align*}
It follows from H\"{o}lder inequality, Young inequality and Lemma \ref{2.1} that
\begin{align}\label{3.3.6}
\nonumber&\left|\int_\Omega\left(F_N(\|u^1\|)(u^1\cdot\nabla)u^1-F_N(\|u^2\|)(u^2\cdot\nabla)u^2\right)\cdot Au\,dx\right|\\
\nonumber\leq&|F_N(\|u^1\|)|\,\|u\|_{L^{\infty}(\Omega)}\,\|\nabla u^1\|_2\,\|Au\|_2+|F_N(\|u^2\|)|\,\|u^2\|_6\,\|\nabla u\|_3\,\|Au\|_2\\
\nonumber&+|F_N(\|u^1\|)-F_N(\|u^2\|)|\,\|u^2\|_6\,\|\nabla u^1\|_3\,\|Au\|_2\\
\nonumber\leq&C|F_N(\|u^1\|)|\,\|u^1\|\,\|u\|^{\frac{1}{2}}\,\|Au\|_2^{\frac{3}{2}}+C|F_N(\|u^2\|)|\,\|u^2\|\,\|u\|^{\frac{1}{2}}\,\|Au\|_2^{\frac{3}{2}}\\
\nonumber&+\frac{C}{N}|F_N(\|u^1\|)|\,|F_N(\|u^2\|)|\,\|u^2\|\,\|u^1\|^{\frac{1}{2}}\,\|Au^1\|_2^{\frac{1}{2}}\,\|u\|\,\|Au\|_2\\
\nonumber\leq&CN\|u\|^{\frac{1}{2}}\,\|Au\|_2^{\frac{3}{2}}+CN^{\frac{1}{2}}\|Au^1\|_2^{\frac{1}{2}}\,\|u\|\,\|Au\|_2\\
\leq& C(N^4+N\|Au^1\|_2)\|u\|^2+\frac{\nu}{2}\|Au(t)\|_2^2.
\end{align}
Therefore, we infer from \eqref{3.3.5}-\eqref{3.3.6} that
\begin{align}\label{3.3.7}
\frac{d}{dt}\|u(t)\|^2+\nu\|Au(t)\|_2^2\leq C(N^4+N\|Au^1\|_2)\|u\|^2.
\end{align}
For any $t\geq t_2,$ multiplying \eqref{3.3.7} by $\bar{t}=t-t_2$ and integrating the resulting inequality over $(t_2,t),$ we obtain
\begin{align*}
\bar{t}\|u(t)\|^2\leq &C\int_{t_2}^t(N^4+N\|Au^1(s)\|_2)(s-t_2)\|u(s)\|^2\,ds+\int_{t_2}^t\|u(s)\|^2\,ds.
\end{align*}
It follows from the classical Gronwall inequality, \eqref{3.3.2} and \eqref{3.3.4} that
\begin{align*}
\bar{t}\|u(t)\|^2\leq &\left(\int_{t_2}^t\|u(s)\|^2\,ds\right)\exp\left\{C\int_{t_2}^t(N^4+N\|Au^1(s)\|_2)\,ds\right\}\\
\leq&\left(\int_0^t\|u(s)\|^2\,ds\right)\exp\left\{CN^4\bar{t}+CN\bar{t}^{\frac{1}{2}}\left(\int_{t_2}^t\|Au^1(s)\|_2^2\,ds\right)^{\frac{1}{2}}\right\}\\
\leq&\frac{1}{\nu}e^{CN^4t}\exp\left\{CN^4\bar{t}+CN\nu^{-\frac{1}{2}}\bar{t}^{\frac{1}{2}}\left(\frac{2}{\nu}\|f\|_2^2\bar{t}+CN^4\rho_1\bar{t}+\rho_1\right)^{\frac{1}{2}}\right\}\|u(0)\|_2^2\\
\leq&\frac{1}{\nu}\exp\left\{CN^4\bar{t}+CN\nu^{-\frac{1}{2}}\bar{t}\left(\frac{2}{\nu}\|f\|_2^2+CN^4\rho_1+\rho_1\right)^{\frac{1}{2}}\right\}\|u(0)\|_2^2
\end{align*}
for any $t\geq t_2+1.$\\
\qed\hfill

Next, we will prove the time regularity of the semigroup $\{S(t)\}_{t\geq0}$ generated by problem \eqref{1.1}. 
\begin{theorem}\label{3.3.8}
Assume that $f\in H.$ Then for any bounded subset $B\subset H,$ there exists a positive function $\varrho_3(s)$ defined on $\mathbb{R}^+$ and some time $t_3=t_3(B)>0$ such that
\begin{align}
\|S(t)u_0-S(\tilde{t})u_0\|_2\leq \varrho_3(|t-\tilde{t}|)|t-\tilde{t}|^{\frac{1}{2}}
\end{align}
for any $t,\tilde{t}\geq t_3$ and any $u_0\in B,$ where $S(t)u_0$ is the solution of problem \eqref{1.1} with initial data $u_0.$
\end{theorem}
\textbf{Proof.} For any bounded subset $B$ of $H,$ we deduce from Theorem \ref{3.2.1} that there exists some time $t_3=t_3(B)$ such that
\begin{align*}
S(t)B\subset B_1=\{u\in V:\|u\|^2\leq \rho_1\}
\end{align*}
for any $t\geq t_3.$

For any $t,\tilde{t}\geq t_3,$ without loss of generality, we assume that $t\geq\tilde{t},$ integrating \eqref{3.2.5} over $(\tilde{t}, t),$ we obtain that for any $u_0\in B,$
\begin{align}\label{3.3.9}
\nonumber\nu\int_{\tilde{t}}^t\|Au(s)\|_2^2\,ds\leq&\frac{2}{\nu}\|f\|_2^2(t-\tilde{t})+CN^4\int_{\tilde{t}}^t\|u(s)\|^2\,ds+\|u(\tilde{t})\|^2\\
\leq&\frac{2}{\nu}\|f\|_2^2(t-\tilde{t})+CN^4\rho_1(t-\tilde{t})+\rho_1.
\end{align}
Thanks to 
\begin{align}\label{3.3.10}
\nonumber\|u_t\|_2\leq&\|f\|_2+\nu\|Au\|_2+|F_N(\|u\|)|\,\|u\|_6\,\|\nabla u\|_3\\
\leq&\|f\|_2+\nu\|Au\|_2+C|F_N(\|u\|)|\,\|u\|\,\|Au\|_2.
\end{align}
Combining Lemma \ref{2.1} with \eqref{3.3.9}-\eqref{3.3.10}, we obtain
\begin{align*}
\int_{\tilde{t}}^t\|u_t(s)\|_2^2\,ds\leq&2\|f\|_2^2(t-\tilde{t})+(2\nu^2+CN^2)\int_{\tilde{t}}^t\|Au(s)\|_2^2\,ds\\
\leq&2\|f\|_2^2(t-\tilde{t})+\frac{1}{\nu}(2\nu^2+CN^2)\left(\frac{2}{\nu}\|f\|_2^2(t-\tilde{t})+CN^4\rho_1(t-\tilde{t})+\rho_1\right)\\
\triangleq:&\varrho_3(|t-\tilde{t}|)^2.
\end{align*}
Therefore, we obtain
\begin{align*}
\|S(t)u_0-S(\tilde{t})u_0\|_2=&\left\|\int_{\tilde{t}}^t u_t(s)\,ds\right\|_2\\
\leq&\int_{\tilde{t}}^t \|u_t(s)\|_2\,ds\\
\leq&\left(\int_{\tilde{t}}^t\|u_t(s)\|_2^2\,ds\right)^{\frac{1}{2}}|t-\tilde{t}|^{\frac{1}{2}}\\
\leq&\varrho_3(|t-\tilde{t}|)|t-\tilde{t}|^{\frac{1}{2}}
\end{align*}
for any $t,\tilde{t}\geq t_3$ and any $u_0\in B.$\\
\qed\hfill

Finally, inspired by the idea in \cite{ba, ea,em, rjc}, we can easily construct the existence of an exponential attractor for problem \eqref{1.1}. Here, we only state it as follows.
\begin{theorem}\label{3.3.11}
Assume that $f\in H.$ Let $\{S(t)\}_{t\geq0}$ be a semigroup generated by problem \eqref{1.1}. Then the semigroup $\{S(t)\}_{t\geq
0}$ possesses an exponential attractor $\mathcal{E}\subset H,$ namely,
\begin{itemize}
\item [(i)] $\mathcal{E}$ is compact and positively invariant with respect to $S(t),$ i.e.,
\begin{align*}
S(t)\mathcal{E}\subset\mathcal{E}
\end{align*}
for any $t\geq0.$
\end{itemize}
\item [(ii)] The fractal dimension $dim_f(\mathcal{E},H)$ of $\mathcal{E}$ is finite.
\item [(iii)] $\mathcal{E}$ attracts exponentially any bounded subset $B$ of $H,$ that is, there exists a positive nondecreasing function $Q$ and a constant $\rho>0$ such that
    \begin{align*}
    dist_H(S(t)B,\mathcal{E})\leq Q(\|B\|_H)e^{-\rho t}
    \end{align*}
    for any $t\geq 0,$ where $dist_H$ denotes the non-symmetric Hausdorff semi-distance between two subsets of $H$ and $\|B\|_H$ stands for the size of $B$ in $H.$ Moveover, both $Q$ and $\rho$ can be explicitly calculated.
\end{theorem}
Thanks to $\mathcal{A}_0\subset\mathcal{E},$ we immediately deduce the following result.
\begin{corollary}
Assume that $f\in H.$ Then the fractal dimension of the global attractor for problem \eqref{1.1} established in Theorem \ref{3.2.9} is finite, i.e.,
\begin{align*}
dim_f(\mathcal{A}_0,H)<+\infty.
\end{align*}
\end{corollary}

\section*{Acknowledgement}
 This work was supported by the National Science Foundation of China Grant (11401459).

\bibliographystyle{elsarticle-template-num}
\bibliography{BIB}

\end{document}